\documentclass[12pt,a4paper]{article}
\usepackage{amsmath,amssymb,amsfonts}
\def\Bbb{\mathbb}

\title{\bf On the equivalence of certain quadratic irrationals}

\author{Kurt Girstmair}

\date{}
\makeatletter
\let\@@maketitle=\maketitle
\def\maketitle{\def\thispagestyle##1{\relax}\@@maketitle}
\makeatother
\newtheorem{theorem}{Theorem}

\newtheorem{lemma}{Lemma}
\newtheorem{corollary}{Corollary}

\def\BE{\begin{equation}}
\def\EE{\end{equation}}
\def\BD{\begin{displaymath}}
\def\ED{\end{displaymath}}
\def\BA{\begin{array}}
\def\EA{\end{array}}
\def\BEA{\begin{eqnarray*}}
\def\EEA{\end{eqnarray*}}
\def\BI{\bibitem}

\def\Z{\Bbb Z}
\def\Q{\Bbb Q}

\def\phi{\varphi}

\def\MB{\mbox}
\def\LD{\ldots}
\def\OV{\overline}

\def\DIV{\,|\,}
\def\NDIV{\, \nmid \,}

\def\MN{\medskip\noindent}
\def\STOP{\hfill$\Box$}

\begin{document}
\maketitle

\begin{abstract}

\noindent
This paper deals with quadratic irrationals of the form $m/q+\sqrt v$ for fixed positive integers $v$ and $q$, $v$ not a square, and varying integers $m$, $(m,q)=1$.
Two numbers $m/q+\sqrt v$, $n/q+\sqrt v$ of this kind  are equivalent (in a
classical sense) if their continued fraction expansions can be written with the same period. We give a necessary and sufficient condition for the equivalence in terms of solutions of Pell's
equation. Moreover, we determine the number of equivalence classes to which these quadratic irrationals belong.
\end{abstract}

%%%%%%%%%%%%%%%%%%%%%%%%%%%%%%%%%%%%%%%%%%%%%
\section*{1. Introduction and results}
%%%%%%%%%%%%%%%%%%%%%%%%%%%%%%%%%%%%%%%%%%%%%

Let $v$ and $q$ be positive integers, $v$ not a square. In this paper we study the equivalence between
numbers
\BD
  x=m/q+\sqrt v,
\ED
where $m$ is an integer, $(m,q)=1$. Thus, $v$ and $q$ are fixed, whereas $m$ may vary.

The equivalence of two numbers $x$, $y$ of this kind means that the (regular) continued fractions of $x$ and $y$
can be written with the same period, say,
\BE
\label{1.0}
 x=[a_0,\LD,a_{j-1},[b_1,\LD,b_k]], y=[c_0,\LD,c_{l-1},[b_1,\LD,b_k]],
\EE
where $[b_1,\LD,b_k]$ is the common period. Here the pre-periods $a_0,\LD,a_{j-1}$ and $c_0,\LD,c_{l-1}$ need not occur. In general, it is more likely that
you find equivalent numbers $x$ and $y$  than inequivalent ones, for example,
\BEA
   x&=&\frac 1{12}+\sqrt 7=[2, [1, 2, 1, 2, 4, 5, 16, 47, 1, 1, 3, 1, 1, 4]],\\
   y&=&\frac 5{12}+\sqrt 7=[3, [16, 47, 1, 1, 3, 1, 1, 4, 1, 2, 1, 2, 4, 5]]\\
    &=&[3,16,47,1,1,3,1,1,4,[1, 2, 1, 2, 4, 5, 16, 47, 1, 1, 3, 1, 1, 4]].
\EEA
We write $x\sim y$ if $x$ and $y$ are equivalent.
It is a classical result of Serret that $x\sim y$ is the same as
\BE
\label{1.1}
  y=\frac{ax+b}{cx+d}, \MB{ where }\left(\begin{array}{cc}
                                               a & b \\
                                               c & d \\
                                          \end{array}
                                   \right)  \in \MB{GL}(2,\Z),
\EE
i.e., $a,b,c,d\in \Z$ and $ad-bc=\pm 1$ (see (see \cite[p. 54]{Pe}, \cite[p. 38]{Borw})).

Our first aim is the following theorem.

\begin{theorem} % Theorem 1 %%%%%%%%%%%%%%%%%%%%%%%%%%%%%%%%%%%%%%%%%%%%%%%%%%%%%%%%%%%%%%%
\label{t1}

Let $x=m/q+\sqrt v$, $y=n/q+\sqrt v$, $(m,q)=(n,q)=1$. Let $q_1=(m-n,q)$. Then $x\sim y$ if, and only if, the
equation
\BE
\label{1.2}
  r^2-c^2v=\pm 1
\EE
has a solution $(r,c)\in\Z^2$ such that $(c,q^2)=qq_1.$

\end{theorem} %%%%%%%%%%%%%%%%%%%%%%%%%%%%%%%%%%%%%%%%%%%%%%%%%%%%%%%%%%%%%%%%%%%%%%%%%%%%%%

Of course, (\ref{1.2}) is known as Pell's equation. Our next question concerns the number of equivalence classes to which our quadratic irrationals belong.
Since $x+1\sim x$, we may restrict ourselves to numbers $x=m/q+\sqrt v$, $y=n/q+\sqrt v$ with $0\le m,n\le q-1$.

\begin{theorem} % Theorem 2 %%%%%%%%%%%%%%%%%%%%%%%%%%%%%%%%%%%%%%%%%%%%%%%%%%%%%%%%%%%%%%%
\label{t2}

Let $q_0$ be the smallest divisor of $q$ such that there is a solution $(r,c)$ of {\rm (\ref{1.2})} with $(c,q^2)=qq_0$. Then the numbers $x=m/q+\sqrt v$, $(m,q)=1$,
$0\le m\le q-1$,
belong to exactly $\phi(q_0)$ equivalence classes, each of which contains $\phi(q)/\phi(q_0)$ elements $x$.

\end{theorem} %%%%%%%%%%%%%%%%%%%%%%%%%%%%%%%%%%%%%%%%%%%%%%%%%%%%%%%%%%%%%%%%%%%%%%%%%%%%%%

\MN
{\em Remarks.} 1. Note that every equivalence class contains many elements different from the numbers $x$ in question. For instance, $1/3+\sqrt 2\sim 9\sqrt 2/2$, the latter not being of the appropriate form
for $v=2$ and $q=3$. Equivalent numbers have the same discriminant (see \cite[p. 41]{Borw}). Since the discriminant of $x=m/q+\sqrt v$, $(m,q)=1$, equals $4q^4v$, $x$ cannot be equivalent to a number
$m'/q'+\sqrt v$, $(m',q')=1$, $q'>0$, $q'\ne q$.

2. The unit group $\Z[\sqrt v]^{\times}$ of the ring $\Z[\sqrt v]$ is isomorphic to $\Z\times \Z/2\Z$ and is generated by a fundamental unit $s+t\sqrt v$ together with $-1$. For $r+c\sqrt v \in\Z[\sqrt v]^{\times}$
we have $q\DIV c$ if, and only if, $r+c\sqrt v \in\Z[q\sqrt v]^{\times}$, the unit group of the subring $\Z[q\sqrt v]$. This group has a finite index $k$ in $\Z[\sqrt v]^{\times}$ (see \cite[p. 296]{Ha}). Accordingly,
$(s+t\sqrt v)^k$ is an element of $\Z[q\sqrt v]^{\times}$.

3. It may happen that $\Z[q\sqrt v]^{\times}$ coincides with $\Z[qq_1\sqrt v]^{\times}$ for some divisor $q_1$ of $q$, $q_1>2$.
In this case $c$ is divisible by $qq_1$ for each $r+c\sqrt v \in \Z[q\sqrt v]^{\times}$. In particular, $qq_1\DIV (c,q^2)$ for all these units.
Let $q_0$ be the smallest divisor of $q$ such that there is a unit $r+c\sqrt v\in\Z[q\sqrt v]^{\times}$ with $(c,q^2)=qq_0$. Then $qq_1\DIV (c,q^2)=qq_0$, and,
accordingly, $q_1\DIV q_0$.
By Theorem \ref{t2}, the numbers $m/q+\sqrt v$, $(m,q)=1$, $0\le m\le q-1$, belong to $\phi(q_0)\ge\phi(q_1)>1$ equivalence classes. In particular, not all of these numbers
are equivalent.

\MN
{\em Example.} Let $v=979$ and $q=12$. The fundamental unit in $\Z[\sqrt v]$ is $s+t\sqrt v$ with $s=360449$ and $t=11520=q^2\cdot 80$.
Hence $\Z[\sqrt v]^{\times}=\Z[q^2\sqrt v]^{\times}$ and for every
$r+c\sqrt v\in \Z[\sqrt v]^{\times}$ we have $q^2\DIV c$.
Accordingly, the number $q_0$ of Theorem \ref{t2} equals $q=12$, and the four numbers
$m/12+\sqrt{979}$, $m\in\{1,5,7,12\}$ belong to four different equivalence classes. Indeed,
\BEA
 \frac 1{12}+\sqrt{979}&=&[31, [2, 1, 2, 5, 2, 3, 6, 1, 4, 62]],\\
 \frac 5{12}+\sqrt{979}&=&[31, [1, 2, 2, 1, 1, 13, 1, 4, 1, 6, 1, 61]].
\EEA
So these numbers have periods of different lengths. The numbers $7/12+\sqrt{979}$ and  $5/12+\sqrt{979}$ have inverse periods, and $11/12+\sqrt{979}$ and $1/12+\sqrt{979}$, too.
In general, we say that $x$ and $y$ have {\em inverse periods} if they can be written as in (\ref{1.0}), the period of $y$ being $[b_k,b_{k-1},\LD,b_1]$, however.

\medskip
Theorem \ref{t1} answers the question whether $x=m/q+\sqrt v$ and $y=n/q+\sqrt v$ have inverse periods. This happens if, and only if, $x\sim y'=n/q-\sqrt v$
(see \cite[p. 77]{Pe}). Since $y'\sim -y'=-n/q+\sqrt v$, we obtain the following corollary to
Theorem \ref{t1}.

\begin{corollary} % Corollary 1 %%%%%%%%%%%%%%%%%%%%%%%%%%%%%%%%%%%%%%%%%%%%%%%%%%%%%%%%%%%%%%%
\label{c1}

Let $x=m/q+\sqrt v$, $y=n/q+\sqrt v$ be as above. Let $q_1'=(m+n,q)$. Then $x$ and $y$ have inverse periods if, and only if, the
equation
{\rm (\ref{1.2})}
has a solution $(r,c)\in\Z^2$ such that $(c,q^2)=qq_1'.$

\end{corollary} %%%%%%%%%%%%%%%%%%%%%%%%%%%%%%%%%%%%%%%%%%%%%%%%%%%%%%%%%%%%%%%%%%%%%%%%%%%%%%

We say that $x$ has a {\em self-inverse}\, period if $x$ can be written with a period $[b_1,\LD,b_k]$ but also with the period $[b_k,\LD,b_1]$ (see \cite[p. 78]{Pe}, \cite{Bu}). From Corollary \ref{c1} we obtain

\begin{corollary} % Corollary 2 %%%%%%%%%%%%%%%%%%%%%%%%%%%%%%%%%%%%%%%%%%%%%%%%%%%%%%%%%%%%%%%
\label{c2}

Let $x=m/q+\sqrt v$ be as above. Put $q_1'=2$ if $q$ is even, and $q_1'=1$, otherwise. Then $x$ has a self-inverse period if, and only if, the
equation
{\rm (\ref{1.2})}
has a solution $(r,c)\in\Z^2$ such that $(c,q^2)=qq_1'.$

\end{corollary} %%%%%%%%%%%%%%%%%%%%%%%%%%%%%%%%%%%%%%%%%%%%%%%%%%%%%%%%%%%%%%%%%%%%%%%%%%%%%%

{\em Remarks.} 1. The reader may consult \cite{Ge}, where quadratic irrationals with self-inverse periods are classified by certain
equivalences.

2. Many examples show the following tendency, for which we have no precise mathematical formulation. Namely, if the numbers $m/q+\sqrt v$, $(m,q)=1$, $0\le m\le q-1$, belong to many equivalence classes,
then their periods are short.
For instance, in the case $v=979$, $q=12$ of the above example we have the largest possible number of equivalence classes, which is $4$.
The corresponding period lengths of $m/q+\sqrt v$ are $10$ or $12$.
If we choose $q=9$ instead, then all numbers $m/q+\sqrt v$ belong to the same equivalence class, and the common period of the $6$ elements $m/q$ has length $78$.

%%%%%%%%%%%%%%%%%%%%%%%%%%%%%%%%%%%%%%%%%%%%%
\section*{2. Proofs}
%%%%%%%%%%%%%%%%%%%%%%%%%%%%%%%%%%%%%%%%%%%%%

{\em Proof of Theorem 1.} Let $x=m/q+\sqrt v$, $y=n/q+\sqrt v$, $(m,q)=(n,q)=1$. First suppose $x\sim y$, i.e., there is a matrix $\left(\begin{array}{cc}
                                               a & b \\
                                               c & d \\
                                          \end{array}
                                   \right)  \in \MB{GL}(2,\Z)$
such that (\ref{1.1}) holds. Then comparison of the coefficients with respect to the $\Q$-basis $(1,\sqrt v)$ of $\Q[\sqrt v]$ shows that (\ref{1.1}) is equivalent to the identities
\BE
\label{2.2}
  a=\frac{c(m+n)}q +d
\EE
and
\BE
\label{2.3}
  b=\frac{-d(m-n)}q-\frac{cm^2}{q^2}+cv.
\EE
Since $b\in\Z$, (\ref{2.3}) implies
\BE
\label{2.4}
  d(m-n)q +cm^2\equiv 0\mod q^2.
\EE
However, $(m,q)=1$, so (\ref{2.4}) requires $c\equiv 0\mod q$. Then $a\in\Z$, by (\ref{2.2}). Let $m-n=q_1m_1$ with $q_1=(m-n,q)$ and $(m_1,q/q_1)=1$. Accordingly, (\ref{2.4}) can be
written
\BE
\label{2.6}
  dm_1qq_1+cm^2\equiv 0\mod q^2.
\EE
If $q_1=q$, this congruence yields $c\equiv 0\equiv qq_1\mod q^2$. If $q_1<q$, we have $qq_1\DIV c$. In this case one easily checks that $(c,q^2)$
must be $qq_1$ (observe that $(d,q)=1$ since $q\DIV c$ and the matrix in question
has determinant $\pm 1$).

Moreover, the condition $ad-bc=\pm 1$ is the same as saying
\BE
\label{2.8}
 d=-\frac{cm}{q}+\sqrt{c^2v\pm 1},
\EE
where we do not fix a sign for the square root.
Hence there is an integer $r$ such that (\ref{1.2}) holds.

Conversely, suppose that $(r,c)$ is a solution of (\ref{1.2}) and $c=qq_1c_1$ for some integer $c_1$ with $(c_1, q/q_1)=1$ (observe $q^2/(qq_1)=q/q_1$).
We define $d$ in such a way that (\ref{2.8}) is satisfied, i.e.,
\BE
\label{2.9}
  d=-\frac{cm}{q}+r=-q_1c_1m+r.
\EE
Then condition (\ref{2.4}) reads
\BD
 (-q_1c_1m+r)(m-n)q+qq_1c_1m^2\equiv 0\mod q^2.
\ED
This congruence is equivalent to the congruence
\BE
\label{2.10}
  c_1mn+rm_1\equiv 0\mod q/q_1.
\EE
Here $(m,q)=(n,q)=(m_1,q/q_1)=1$. Observe that $(c_1,q/q_1)=1$. By (\ref{1.2}), we have $(r,q)=1$, since $q\DIV c$.
Of course, it may happen that our pair $(r,c)$ does not satisfy (\ref{2.10}). In this case we consider $r'+c'\sqrt v=(r+c\sqrt v)^k$ for some positive integer $k$ prime to $q$.
Since $q\DIV c$, we obtain
\BD
 r'+c'\sqrt v\equiv r^k+kr^{k-1}c\sqrt v\mod q^2,
\ED
a congruence mod $\Z[\sqrt v]q^2$. It is easy to see that this congruence implies the congruences
\BE
\label{2.12}
r'\equiv r^k\mod q^2,\enspace c'\equiv kr^{k-1}c \mod q^2,
\EE
which are congruences mod $\Z q^2$. In particular, $(c',q^2)=(kc,q^2)=(c,q^2)$, since $(k,q)=1$.
The second of the congruences (\ref{2.12}) shows that we may write $c'=qq_1c_1'$ with $c_1'\equiv kr^{k-1}c_1\mod q/q_1$.
The congruence (\ref{2.10}), for $r'$ and $c'$ instead of $r$ and $c$, reads
\BE
\label{2.13}
  c_1'nm+r'm_1\equiv 0\mod q/q_1,
\EE
or
\BD
kr^{k-1}c_1nm+r^km_1\equiv 0\mod q/q_1.
\ED
Because $(r,q)=1$, this is equivalent to
\BE
\label{2.14}
   kc_1nm+rm_1\equiv 0\mod q/q_1.
\EE
Observe $(c_1,q/q_1)=(m_1,q/q_1)=(m,q)=(n,q)=1$ and $q/q_1\DIV q$. Therefore, the number $k$ (prime to $q$) can be chosen such that (\ref{2.14}) holds.
Then (\ref{2.13}) holds, and, thus, the congruence (\ref{2.4}).

We define $d$ by (\ref{2.9}) with $r'$, $c'$ instead of $r$, $c$.
Finally, we define $a$, $b$ by (\ref{2.2}) and (\ref{2.3}) with $c'$ instead of $c$. Then $a$, $b$ are integers, $ad-bc'=\pm 1$, and (\ref{1.1}) also holds.
\STOP

\MN
{\em Example.} Let $v=7$, $q=12$, $m=1$, $n=5$. Hence $x=1/12+\sqrt{7}$, $y=5/12+\sqrt{7}$, $m-n=-4=q_1m_1$ with $q_1=4,m_1=-1$ (see the example at the beginning of this paper).
The fundamental unit of $\Z[\sqrt{7}]$ is $s+t\sqrt{7}=8+3\sqrt{7}$. Since $(s+t\sqrt 7)^2=127+48\sqrt 7$, we put $r=127$, $c=48$. In particular, $(c,q^2)=qq_1=48$.
So $x$ is equivalent to $y$ (what we already know).
We have $c_1=1$ and $q/q_1=3$. However, the congruence (\ref{2.10}) does not hold, but (\ref{2.14}) is true with $k=5$, $(k,q)=1$.
This choice, however, leads to rather large numbers. But we may also choose $k=2$ (which is not prime to $q$).
Indeed, define $r'$, $c'$ by
\BD
 r'+c'\sqrt{7}=(r+c\sqrt{7})^{2}=32257+12192\sqrt{7}.
\ED
Then $(c',q^2)=(12192,144)=48=qq_1$, as required.
As in the proof of Theorem \ref{t1} we obtain
$d=31241$, $a=37337$, \MB{ and } $b=95673$.
In this way $ad-bc'=1$ and $(ax+b)/(c'x+d)=y$.

%see "Cfrac_22_14.mw"

\bigskip
The proof of Theorem \ref{t2} requires the following lemmas.

\begin{lemma} % Lemma 2 %%%%%%%%%%%%%%%%%%%%%%%%%%%%%%%%%%%%%%%%
\label{l2}

Let $(r_1,c_1)$ and $(r_2,c_2)$ be solutions of {\rm (\ref{1.2})} such that $(c_1,q^2)=qq_1$ and $(c_2,q^2)$ $=qq_2$ for divisors $q_1$, $q_2$ of $q$.
Let $q'=(q_1,q_2)$. Then there is a solution $(r',c')$ of {\rm (\ref{1.2})} such that $(c',q^2)=q'$.

\end{lemma} %%%%%%%%%%%%%%%%%%%%%%%%%%%%%%%%%%%%%%%%%%%%%%%%%%%

\MN
{\em Proof.} Let $j_1,j_2$ be positive integers. We have
\BD
  (r_i+c_i\sqrt v)^{j_i}\equiv r_i^{ji}+j_ir_i^{j_i-1}c_i\sqrt v \mod q^2,
\ED
$i=1,2$. We define $r'+c'\sqrt v$ by
\BD
   r'+c'\sqrt v=(r_1+c_1\sqrt v)^{j_1}(r_2+c_2\sqrt v)^{j_2}.
\ED
Since
\BD
 (r_1+c_1\sqrt v)^{j_1}(r_2+c_2\sqrt v)^{j_2}\equiv r_1^{j_1}r_2^{j_2}+(j_1r_1^{j_1-1}r_2^{j_2}c_1+j_2r_1^{j_1}r_2^{j_2-1}c_2)\sqrt v\mod q^2,
\ED
we obtain
\BE
\label{2.16}
  c'\equiv r_1^{j_1-1}r_2^{j_2-1}(j_1r_2c_1+j_2r_1c_2)\mod q^2.
\EE
Observe that $(r_1,q)=(r_2,q)=1$, since $q\DIV c_1,c_2$.
We consider the ideal $\Z r_2c_1+\Z q^2$ in $\Z$. We have $(r_2c_1,q^2)=qq_1$. This can be written as
\BD
  \Z r_2c_1+\Z q^2=\Z qq_1.
\ED
In the same way, we obtain
\BD
  \Z r_1c_2+\Z q^2=\Z qq_2.
\ED
However, $(qq_1,qq_2)=qq'$, and so $\Z qq_1+\Z qq_2=\Z qq'$. This yields
\BD
  \Z r_2c_1+\Z r_1c_2 +\Z q^2=\Z qq'.
\ED
Accordingly, there are integers $k_1$, $k_2$ such that
\BE
\label{2.17}
 k_1r_2c_1+k_2r_1c_2 \equiv qq'\mod q^2.
\EE
We put $j_1=k_1+lq^2$ and $j_2=k_2+lq^2$, where $l$ is chosen such that both $j_1,j_2$ are positive.
Then (\ref{2.16}) and (\ref{2.17}) show that $c'$ satisfies $(c',q^2)=qq'$.
\STOP

\begin{lemma} % Lemma 3 %%%%%%%%%%%%%%%%%%%%%%%%%%%%%%%%%%%%%%%%
\label{l3}

Let $q_1$ divide $q$.
Let $(r,c)$ be a solution of {\rm (\ref{1.2})} such that $(c,q^2)=qq_1$.
Let $p$ be a prime dividing $q/q_1$. Define $(r',c')$ by $r'+c'\sqrt v=(r+c\sqrt v)^p$.
Then $(r',c')$ is a solution of {\rm (\ref{1.2})} such that $(c',q^2)=qq_1p.$

\end{lemma} %%%%%%%%%%%%%%%%%%%%%%%%%%%%%%%%%%%%%%%%%%%%%%%%%%%

\MN
{\em Proof.} Let $l$ be a prime number. By $v_l(k)$ we denote the $l$-exponent of the integer $k$, i.e., $l^{v_l(k)}\DIV k$, $l^{v_l(k)+1}\NDIV k$.
If $l$ divides $q$, we have $v_l(c)=v_l(qq_1)$, since $qq_1\DIV q^2$ and $(c,q^2)=qq_1$.

First let $p=2$. Then
$c'=2rc$ and $v_2(2rc)=v_2(c)+1$, since $v_2(r)=0$, by (\ref{1.2}).
If $l$ is a prime divisor of $q$ different from $2$, we see $v_l(2rc)=v_l(c)$.
Hence $(c',q^2)=qq_1p$.

If $p\ge 3$ we have
\BE
\label{2.18}
  r'+c'\sqrt v\equiv r^p+pr^{p-1}c\sqrt v \mod c^2p,
\EE
since $c^p\equiv 0\mod c^2p$ (recall $p\DIV c$). Thus, $v_p(c')=v_p(pr^{p-1}c)=v_p(c)+1$, because $v_p(c^2p)>v_p(c)+1$. For a prime divisor $l$ of $q$ different from $p$ we have
$v_l(c^2p)=2v_l(c)>v_l(c) (\ge 1)$. From (\ref{2.18}) we obtain $v_l(c')=v_l(pr^{p-1}c)=v_l(c)$.
\STOP

\MN
{\em Proof of Theorem \ref{t2}.}
Let $q_0$ be the smallest divisor of $q$ such that there is a solution $(r,c)$ of (\ref{1.2}) with $(c,q^2)=qq_0$.

If $q$ is even,  then $m-n$ is even for all $m$, $n$ with $(m,q)=(n,q)=1$. Hence $(m-n,q)$ is even. Accordingly, $q_0$ cannot have the
form $q_0=(m-n,q)$ if $q_0$ is odd. Suppose that this holds. Then we replace $q_0$ by $2q_0$.
Since $\phi(2q_0)=\phi(q_0)$, this does not change the assertion of Theorem \ref{t2}.
Moreover, by Lemma \ref{l3}, we have a solution $(r',c')$ of (\ref{1.2}) such that $(c',q^2)=2qq_0$.

Accordingly, we may assume that $q_0$ is even if $q$ is even and
suppose that $(r,c)$ is a solution of (\ref{1.2}) with $(c,q^2)=qq_0$.

Let $y=n/q+\sqrt v, (n,q)=1$. We show that the sets
\BD
  X_1=\{m/q+\sqrt v: (m,q)=1, m/q+\sqrt v\sim y\}
\ED
and
\BD
 X_2= \{m/q+\sqrt v: (m,q)=1, q_0\DIV m-n\}
\ED
coincide. Indeed, if $m/q+\sqrt v$ is in $X_1$ and $(m-n,q)=q_1$, then there is a solution $(r',c')$ of (\ref{1.2}) such that $(c',q^2)=qq_1$. On the other hand,
we have a solution $(r,c)$ of (\ref{1.2}) such that $(c,q^2)=qq_0$. By Lemma \ref{l2}, there is a solution $(r'',c'')$ such that $(c'',q^2)=q(q_0,q_1)$.  If $(q_0,q_1)\ne q_0$, then
$(q_0,q_1)<q_0$, which contradicts the minimality of $q_0$. Accordingly, $(q_0,q_1)=q_0$ and $q_0\DIV q_1=(m-n,q)$. In particular, $q_0$ divides $m-n$ and $m/q+\sqrt v\in X_2$.

Conversely, if $m/q+\sqrt v\in X_2$, then $(m-n,q)=q_0k$ for some positive integer $k$. By Lemma \ref{l3}, there exists a solution $(r',c')$ of (\ref{1.2}) such that $(c',q^2)=q_0k$.
This implies $m/q+\sqrt v\sim y$ and $m/q+\sqrt v\in X_1$.

We consider $X=\{m/q+\sqrt v: (m,q)=1, 0\le m\le q-1\}$ and $X_1'=X_1\cap X$.
%The set $X_1'$ consist of all numbers $x\in X$ that are equivalent to $y$.
Now $m/q+\sqrt v\in X$ lies in $X_1'$ if, and only if, $q_0\DIV m-n$, i.e., the canonical surjection
\BD
   \pi:(\Z/q\Z)^{\times}\to (\Z/q_0\Z)^{\times}: \OV k\mapsto \OV k
\ED
maps $\OV m$ onto $\OV n$. Thereby, $|X_1'|=|\pi^{-1}(\OV n)|=\phi(q)/\phi(q_0)$.
Hence there are exactly $\phi(q)/\phi(q_0)$ elements of $X$ that are equivalent to $y$. This, however, implies that there must be exactly $\phi(q_0)$ equivalence classes whose intersections with $X$ are not empty.
\STOP

\MN
\centerline{\bf Dataset availability}

\MN
The datasets generated during and analysed during the current study are available from the author on reasonable request.

%%%%%%%%%%%%%%%%%%%%%%%%%%%%%%%%%%%%%%%%%%%%%%%%%%%%%
%%%%%%%%%%%%%%%%%%%%%%%%%%%%%%%%%%%%%%%%%%%%%%%%%%%%%%%%%%%%%%%%%%%%%%%%%%

\MN
Kurt Girstmair\\
Institut f\"ur Mathematik \\
Universit\"at Innsbruck   \\
Technikerstr. 13/7        \\
A-6020 Innsbruck, Austria \\
Kurt.Girstmair@uibk.ac.at

\end{document}